\documentclass{elsart}

\journal{}

\usepackage{amsmath}
\usepackage{amssymb}

\begin{document}

\newcommand {\vt}{\vartheta}
\newcommand {\po}[2]{({#1};{#2})_\infty}
\newtheorem {Theorem}{Theorem}
\newtheorem {Corollary}{Corollary}
\newtheorem {Definition}{Definition}

\begin{frontmatter}

 \title{Several special values of Jacobi theta functions}
 \author{Istv\'an Mez\H{o}\thanksref{a}\thanksref{b}}
 \thanks[a]{Present address: Departamento de Matem\'atica\\Escuela Polit\'ecnica Nacional\\Ladr\'on de Guevara E11-253\\Quito, Ecuador}
 \thanks[b]{This scientific work was financed by Proyecto Prometeo de la Secretar\'ia Nacional de Ciencia, Tecnolog\'ia e Innovaci\'on (Ecuador).}
 \ead{istvan.mezo@epn.edu.ec}
 \ead[url]{http://www.inf.unideb.hu/valseg/dolgozok/mezoistvan/mezoistvan.html}

\begin{abstract}Using the duplication formulas of the elliptic trigonometric functions of Gosper, we deduce some new special values for the first two Jacobi theta functions. At the end of the paper, we show how is it possible to extend our arguments and deduce a wide variety of additional special values for the Jacobi thetas. In addition, an identity is revealed between these functions and the Weber modular function.
\end{abstract}

\end{frontmatter}

\section{Introduction}

By using the relatively new definition of a kind of $q$-trigonometric function due to R. W. Gosper, we deduce a number of closed form evaluations of the classical Jacobi theta functions $\vt_1(z,q)$ and $\vt_2(z,q)$ and some of their quotients. We mention a few results preliminarily:
\begin{align*}
i\vt_1\left(-\frac{i\pi}{4},e^{-2\pi}\right)&=2^{5/8}\frac{e^{\pi /32}}{\pi ^{3/4}}\Gamma \left(\frac{5}{4}\right) \sqrt[16]{-68+323 \sqrt{2}-56 \sqrt{68-14 \sqrt{2}}},\\
i\vt_1\left(-\frac{i\pi}{2},e^{-4\pi}\right)&=\frac{e^{\pi /16}}{\pi ^{3/4}}\Gamma \left(\frac{5}{4}\right)\sqrt[4]{2+6 \sqrt{2}-4 \sqrt{2 \left(2+\sqrt{2}\right)}} .
\end{align*}
Our main tool will be a relation between the $\vt_1$ function and the $q$-sine function of Gosper, and a duplication formula for the latter. Both results due to the author \cite{MezoPAMS}.

Later, we reveal a combinatorial relation between two elliptic sums, one of them was investigated by Ramanujan. This permits us to connect ratios of the second Jacobi theta function $\vt_2$ to the Weber modular function and to the Dedekind eta. Using the known special values of the latters, we will prove among others that
\[\frac{\vt_2(e^{-\pi/4})}{\vt_2(e^{-\pi/2})}=\left(8+6\sqrt2\right)^{1/8},\]
\[\frac{\vt_2(e^{-\pi/3})}{\vt_2(e^{-2\pi/3})}=\]
\[\frac{\left(594100+420099 \sqrt{2}\right)^{1/12} \left(1+\sqrt{3}+\sqrt{2} 3^{3/4}\right)^{4/3} \left(-2-\sqrt{3}+\sqrt{9+6 \sqrt{3}}\right)}{2^{13/24} \left(-6+5 \sqrt{2}+6\ 3^{1/4}+9 \sqrt{2} 3^{1/4}-4\ 3^{3/4}+\sqrt{2} 3^{3/4}+7 \sqrt{6}\right)^{2/3}}.\]

\section{Preliminaries}

\subsection{The Jacobi theta functions}

The first and second Jacobi theta functions are defined via the next doubly infinite sums \cite{ArmitageEberlein}:
\begin{eqnarray}
\vartheta_1(z,q)&=&\sum_{n=-\infty}^\infty(-1)^{n-\frac12}q^{\left(n+\frac12\right)^2}e^{(2n+1)iz},\label{th1dsum}\\
\vartheta_2(z,q)&=&\sum_{n=-\infty}^\infty q^{\left(n+\frac12\right)^2}e^{(2n+1)iz}\label{th2dsum},
\end{eqnarray}
where $i=\sqrt{-1}$, $z\in\mathbb{C}$ and $|q|<1$.

There are representations of the above functions only with simply infinite (unilateral) sums:
\begin{eqnarray}
\vartheta_1(z,q)&=&2\sum_{n=0}^\infty(-1)^nq^{\left(n+\frac12\right)^2}\sin[(2n+1)z],\label{th1ssum}\\
\vartheta_2(z,q)&=&2\sum_{n=0}^\infty q^{\left(n+\frac12\right)^2}\cos[(2n+1)z]\label{th2ssum}.
\end{eqnarray}
We use the abbreviation $\vartheta_i(0,q)=\vartheta_i(q)$ ($i=1,2$).

\subsection{The Gosper $q$-trigonometric functions}

Gosper defined his $q$-trigonometric functions as follows \cite{Gosper}:
\begin{eqnarray}
\sin_q(\pi z)&:=&q^{(z-1/2)^2}\frac{(q^{2z};q^2)_\infty(q^{2-2z};q^2)_\infty}{(q;q^2)_\infty^2}\quad(0<q<1),\label{qsindef}\\
\cos_q(\pi z)&:=&q^{z^2}\frac{(q^{1-2z};q^2)_\infty(q^{2z+1};q^2)_\infty}{(q;q^2)_\infty^2}\quad(0<q<1).\label{qcosdef}
\end{eqnarray}
It can be seen that $\cos_q(z)=\sin_q\left(\frac{\pi}{2}-z\right)$. Moreover
\begin{equation}
\sin_q(z)=\frac{\vartheta_1(z,p)}{\vartheta_2(p)},\quad \cos_q( z)=\frac{\vartheta_2(z,p)}{\vartheta_2(p)}.\label{qsinth}
\end{equation}
Here $p$ is implicitly defined by the equation
\begin{equation}
\ln p\cdot\ln q=\pi^2.\label{pdef}
\end{equation}

Based on computer experience, Gosper conjectured \cite{Gosper} that the next theorem is true.

\begin{Theorem}\label{thmsindup}For all $0<q<1$ and $z\in\mathbb{R}$, we have
\begin{equation}
\sin_q(2z)=C(q)\sin_{q^2}(z)\cos_{q^2}(z),\label{qsineaddM}
\end{equation}
where
\begin{equation}
C(q)=q^{-\frac14}\frac{(q^2;q^4)_\infty^4}{(q;q^2)_\infty^2}=q^{-\frac14}\frac{\phi^6(q^2)}{\phi^2(q)\phi^4(q^4)}.\label{cdef}
\end{equation}
\end{Theorem}

The present author offered a proof \cite{MezoPAMS}, hence now the above is really a theorem. However, this is not the case with the other conjecture of Gosper with respect to the $q$-cosine. He stated that
\begin{equation}
\cos_q(2z)=\cos_{q^2}^2(z)-\sin_{q^2}^2(z).\label{qcosineaddM}
\end{equation}
We could not prove this statement, but the numerical evidences are so strong that we have no any doubt on its validity.

The names $q$-sine and $q$-cosine are supported by the fact that
\[\lim_{q\to1-}\sin_q(z)=\sin(z),\quad\mbox{and}\quad\lim_{q\to1-}\cos_q(z)=\cos(z).\]

See \cite{Gosper,GosperComp,GIZ,MezoPAMS} for more on th Gosper $q$-trigonometric functions.

\subsection{The Euler modular function and the $q$-Pochhammer symbols}

Let $(x;q)_n:=(1-x)(1-qx)\cdots(1-q^{n-1}x)$ and $(x;q)_\infty:=\lim_{n\rightarrow\infty}(x;q)_n$, if the limit exists. In what follows, $q$ will always be a parameter with $|q|<1$. This guarantees that the infinite product
\[\phi(q):=(q;q)_\infty=\prod_{n=1}^\infty(1-q^n)\]
converge. This function is called the Euler (modular) function. S. Ramanujan extensively studied this function as a special case of his theta function \cite{RNB4,RNB3,RNB5}. What we will widely use, are the special values
\begin{align}
\phi\left(e^{-\pi/2}\right)&=\frac{\sqrt[3]{\sqrt{2}-1} \sqrt[24]{4+3 \sqrt{2}} e^{\pi /48} \Gamma \left(\frac{1}{4}\right)}{2^{5/6} \pi ^{3/4}},\label{PhiPp2}\\
\phi\left(e^{-\pi}\right)&=\frac{e^{\pi /24} \Gamma \left(\frac{1}{4}\right)}{2^{7/8} \pi ^{3/4}},\\
\phi\left(e^{-2\pi}\right)&=\frac{e^{\pi /12} \Gamma \left(\frac{1}{4}\right)}{2 \pi ^{3/4}},\\
\phi\left(e^{-4\pi}\right)&=\frac{e^{\pi /6} \Gamma \left(\frac{1}{4}\right)}{2^{11/8} \pi ^{3/4}},\\
\phi\left(e^{-8\pi}\right)&=\frac{\sqrt[4]{\sqrt{2}-1} e^{\pi /3} \Gamma \left(\frac{1}{4}\right)}{2^{29/16} \pi ^{3/4}}.\label{PhiP8}
\end{align}
These relations -- except the first one -- can be found in \cite[p. 326]{RNB5}. The first special value will be deduced later in this paper.

For further reference we note the simple fact
\begin{equation}
\po{q}{q^2}=\frac{\phi(q)}{\phi(q^2)}.\label{qq2toq}
\end{equation}

\subsection{The Jackson $q$-gamma function}

The Jackson $q$-gamma function \cite{Andrews,Gasper} is defined as
\begin{equation}
\Gamma_q(z):=\frac{\phi(q)}{(q^z;q)_\infty}(1-q)^{1-z}\quad(0<q<1).\label{gdammadef}
\end{equation}
This really generalizes the ordinary Euler gamma function, since
\[\Gamma_q(z)\to\Gamma(z)\quad\mbox{if } q\to1-.\]

See more on the $q$-gamma function in \cite[p. 460-]{SC}.

What is important for us is that there exists a simple connection between the $q$-sine function and $\Gamma_q$ \cite{MezoPAMS}:
\begin{equation}
\sin_q(\pi z)=q^{\frac14}\Gamma_{q^2}^2\left(\frac12\right)\frac{(q^2)^{\binom{z}{2}}}{\Gamma_{q^2}(z)\Gamma_{q^2}(1-z)}.\label{fsinqrep}
\end{equation}
This generalizes Euler's reflection formula
\[\sin(\pi z)=\frac{\pi}{\Gamma(z)\Gamma(1-z)},\]
where $\Gamma$ is the Euler gamma function. On an approximative version of \eqref{fsinqrep}, we refer to \cite[Section 7.]{EgM}.

In the next sections we start to build up the ingredients to find special values for $\vt_1$. 

Our argument is the following: by using the addition formula in Theorem \ref{thmsindup}, with some algebra we find several special values for the $\sin_q$ function expressed by the Euler function and then we put these together into the next theorem below to arrive at our goal.

\section{Another formula connecting the first Jacobi theta function and $\sin_q$}

By \eqref{qsinth} we see that the $\sin_q$ Gosper function and $\vt_1$ are intimately connected. We give a new relation now.

\begin{Theorem}\label{theta1sinq}
\begin{equation}
i\vartheta_1(iz\ln q,q)=q^{-z^2}\frac{\phi^2(q)}{\phi(q^2)}\sin_q(\pi z).\label{mid}
\end{equation}
\end{Theorem}

\textit{Proof.} The definition of $\Gamma_q(z)$ gives that
\[\Gamma_{q^2}\left(\frac12\log_q(qy)\right)\Gamma_{q^2}\left(\frac12\log_q(q/y)\right)=\]
\[\frac{(q^2;q^2)_\infty^2(1-q^2)}{(qy;q^2)_\infty(q/y;q^2)_\infty}=\frac{(q^2;q^2)_\infty^3(1-q^2)}{(q^2;q^2)_\infty(qy;q^2)_\infty(q/y;q^2)_\infty}.\]
The denominator can be rewritten via Jacobi's triple product identity \cite[p. 15]{Gasper}:
\[(q^2;q^2)_\infty(qy;q^2)_\infty(q/y;q^2)_\infty=\sum_{n=-\infty}^\infty(-1)^nq^{n^2}y^n.\]
Thus we get that
\[\frac{(q^2;q^2)_\infty^3(1-q^2)}{\Gamma_{q^2}\left(\frac12\log_q(qy)\right)\Gamma_{q^2}\left(\frac12\log_q(q/y)\right)}=\sum_{n=-\infty}^\infty(-1)^nq^{n^2}y^n.\]
Realize that if $z=\frac12\log_q(qy)$, then $1-z=\frac12\log_q(q/y)$ and $y=q^{2z-1}$. Hence
\begin{equation}
\frac{(q^2;q^2)_\infty^3(1-q^2)}{\Gamma_{q^2}(z)\Gamma_{q^2}(1-z)}=\sum_{n=-\infty}^\infty(-1)^nq^{n^2}q^{n(2z-1)}.\label{gamma_series}
\end{equation}
The special value
\begin{equation}
\Gamma_{q^2}\left(\frac12\right)=\frac{(q^2;q^2)_\infty^2}{(q;q)_\infty}\sqrt{1-q^2}\label{gspec}
\end{equation}
and equation \eqref{fsinqrep} yield the remarkable summation formula for $\sin_q$:
\begin{equation}
\sin_q(\pi z)=q^{\left(z-\frac12\right)^2}\frac{(q^2;q^2)_\infty}{(q;q)_\infty^2}\sum_{n=-\infty}^\infty(-1)^nq^{n^2}q^{n(2z-1)}.\label{tmp1}
\end{equation}
The sum on the right can be easily transformed to the theta function $\vartheta_1$:
\begin{equation}
\sum_{n=-\infty}^\infty(-1)^nq^{n^2}q^{n(2z-1)}=q^{-\frac14}\sum_{n=-\infty}^\infty(-1)^nq^{\left(n+\frac12\right)^2}q^{2n(z-1)}.\label{tmp2}
\end{equation}
On the other hand, definition \eqref{th1dsum} gives that
\[ie^{-iy}\vartheta_1(y,q)=\sum_{n=-\infty}^\infty(-1)^nq^{\left(n+\frac12\right)^2}e^{2niy}.\]
If we choose $y=i(1-z)\ln q$, then
\begin{equation}
ie^{-iy}\vartheta_1(y,q)=iq^{1-z}\vartheta_1(i(1-z)\ln q,q)=\sum_{n=-\infty}^\infty(-1)^nq^{\left(n+\frac12\right)^2}q^{2n(z-1)}.\label{theta1series}
\end{equation}
This and \eqref{tmp1}-\eqref{tmp2} gives that
\[\sin_q(\pi z)=q^{\left(z-\frac12\right)^2}\frac{(q^2;q^2)_\infty}{(q;q)_\infty^2}q^{-\frac14}iq^{1-z}\vartheta_1(i(1-z)\ln q,q).\]
Then a straightforward rearrangement and substitution $z=1-z$ give our desired formula \eqref{mid}.
\hfill$\Box$

\section{Special values for $\sin_q$}

In what follows we derive an expression for $\sin_q(\pi/4)$ and $\sin_q(\pi/8)$.

\begin{Theorem}\label{sinqspecval}We have
\begin{eqnarray}
\sin_{q^2}\left(\frac{\pi}{4}\right)&=&q^{\frac18}\frac{(q;q^2)_\infty}{(q^2;q^4)_\infty^2}=q^{\frac18}\frac{\phi(q)}{\phi(q^2)}=\frac{1}{\sqrt{C(q)}},\label{pp4}\\
\sin_{q^2}\left(\frac{\pi}{8}\right)&=&\sqrt{\sin_q\left(\frac{\pi}{4}\right)\frac{\sqrt{1+4/C^2(q)}-1}{2}}.\label{pp8}
\end{eqnarray}
Here $C(q)$ is defined in \eqref{cdef}.
\end{Theorem}
Latter identity shows that we can not wait for a simple closed form expression for $\sin_q(\pi/2^n$) (and $\sin_q(k\pi/2^n)$ in general). However, \eqref{pp8} tends to $\sqrt{\frac{\sqrt{2}}{4}(\sqrt{2}-1)}$, which is $\sin(\pi/8)$ (this is true, because $C(q)\to2$ if $q\to1-$, see the duplication formula in Theorem \ref{thmsindup}).

\textit{Proof.} In order to prove \eqref{pp4}, we put $z=\pi/4$ in \eqref{qsineaddM} and \eqref{qcosineaddM}. Since $\sin_q(\pi/2)=1$ and $\cos_q(\pi/4)=0$ (see the definitions), we easily have that
\[\sin_{q^2}\left(\frac{\pi}{4}\right)=\cos_{q^2}\left(\frac{\pi}{4}\right)=\frac{1}{\sqrt{C(q)}}.\]
According to the definition of $C(q)$, a simple algebraic manipulation proves \eqref{pp4}.

Now we deal with \eqref{pp8}. We apply our duplication formulas for $z=\pi/8$:
\begin{eqnarray*}
\sin_q\left(\frac{2\pi}{8}\right)&=&C(q)\sin_{q^2}\left(\frac{\pi}{8}\right)\cos_{q^2}\left(\frac{\pi}{8}\right),\\
\cos_q\left(\frac{2\pi}{8}\right)&=&\cos_{q^2}^2\left(\frac{\pi}{8}\right)-\sin_{q^2}^2\left(\frac{\pi}{8}\right).\\
\end{eqnarray*}
Since $\cos_q(z)=\sin_q(\pi/2-z)$, we get that the left hand sides coindice, and $\cos_{q^2}(\pi/8)=\sin_{q^2}(3\pi/8)$. Thus
\begin{equation}
C(q)\sin_{q^2}\left(\frac{\pi}{8}\right)\sin_{q^2}\left(\frac{3\pi}{8}\right)=\sin_{q^2}^2\left(\frac{3\pi}{8}\right)-\sin_{q^2}^2\left(\frac{\pi}{8}\right).\label{tmp4}
\end{equation}

Since $\sin_q(z)=\sin_q(\pi-z)$ ($0\le z\le\pi$), the duplication formula for $\sin_q$ with $z=3\pi/8$ easily gives that
\[\sin_{q^2}\left(\frac{3\pi}{8}\right)=\frac{\sin_q\left(\frac{\pi}{4}\right)}{C(q)\sin_{q^2}\left(\frac{\pi}{8}\right)}.\]
Substituting this into \eqref{tmp4}, we win
\[\sin_q\left(\frac{\pi}{4}\right)=\frac{\sin_q^2\left(\frac{\pi}{4}\right)}{C(q)^2\sin_{q^2}^2\left(\frac{\pi}{8}\right)}-\sin_{q^2}^2\left(\frac{\pi}{8}\right).\]
And this is just a quadratic equation in the variable $x=\sin_{q^2}^2\left(\frac{\pi}{8}\right)$. Since the value $\sin_{q^2}\left(\frac{\pi}{8}\right)$ is surely real and positive, the only one solution will be exactly \eqref{pp8}.
\hfill$\Box$

As an interesting corollary, \eqref{pp4} and \eqref{fsinqrep} (together with \eqref{gspec}) yield the curiuous product formula
\[\Gamma_{q^2}\left(\frac14\right)\Gamma_{q^2}\left(\frac34\right)=(1-q^2)\frac{(q^2;q^2)_\infty^4(q;q^2)_\infty^2}{(q;q)_\infty^2(\sqrt{q};q)_\infty}.\]
Or, applying the Euler $\phi$ function and \eqref{qq2toq},
\[\Gamma_{q^2}\left(\frac14\right)\Gamma_{q^2}\left(\frac34\right)=(1-q^2)\frac{\phi^2(q^2)\phi(q)}{\phi(\sqrt{q})}.\]
This is the $q$-analogue of the classical product
\[\Gamma\left(\frac14\right)\Gamma\left(\frac34\right)=\pi\sqrt{2}.\]

Putting the special values of the $\phi$ function into the above formula, we get the next special values.

\begin{Theorem}For the Jackson $q$-gamma function we have
\begin{align*}
\Gamma_{e^{-2\pi}}\left(\frac14\right)\Gamma_{e^{-2\pi}}\left(\frac34\right)&=\frac{e^{-29 \pi /16} \left(e^{2 \pi }-1\right) \Gamma \left(\frac{1}{4}\right)^2}{4\pi ^{3/2} \sqrt[3]{\sqrt{2}-1} \sqrt[24]{8+6 \sqrt{2}}},\\
\Gamma_{e^{-4\pi}}\left(\frac14\right)\Gamma_{e^{-4\pi}}\left(\frac34\right)&=\frac{e^{-29 \pi /8} \left(e^{4 \pi }-1\right) \Gamma \left(\frac{1}{4}\right)^2}{2^{23/8} \pi ^{3/2}},\\
\Gamma_{e^{-8\pi}}\left(\frac14\right)\Gamma_{e^{-8\pi}}\left(\frac34\right)&=\frac{\sqrt{\sqrt{2}-1} e^{-29 \pi /4} \left(e^{8 \pi }-1\right) \Gamma \left(\frac{1}{4}\right)^2}{16 \pi ^{3/2}}.
\end{align*}
\end{Theorem}

By our knowledge, these values are entirely new.

We finish this section with a short remark on the individual values $\Gamma_q\left(\frac12\right)$. However one can prove the next theorem without using any new results of this paper, we had never met explicitly with these values. Hence, by further references we add the following theorem to the paper.

\begin{Theorem}The next special values hold
\begin{align*}
\Gamma_{e^{-\pi}}\left(\frac12\right)&=\frac{e^{-7 \pi /16} \sqrt{e^{\pi }-1} \Gamma \left(\frac{1}{4}\right)}{2^{11/12}\pi ^{3/4} \sqrt[3]{\sqrt{2}-1} \sqrt[24]{4+3 \sqrt{2}}},\\
\Gamma_{e^{-2\pi}}\left(\frac12\right)&=\frac{e^{-7 \pi /8} \sqrt{e^{2 \pi }-1} \Gamma \left(\frac{1}{4}\right)}{2 \sqrt[8]{2} \pi ^{3/4}},\\
\Gamma_{e^{-4\pi}}\left(\frac12\right)&=\frac{e^{-7 \pi /4} \sqrt{e^{4 \pi }-1} \Gamma \left(\frac{1}{4}\right)}{2^{7/4} \pi^{3/4}},\\
\Gamma_{e^{-8\pi}}\left(\frac12\right)&=\frac{e^{-7 \pi /2} \sqrt{\left(\sqrt{2}-1\right) \left(e^{8 \pi }-1\right)} \Gamma \left(\frac{1}{4}\right)}{4 \sqrt[4]{2} \pi ^{3/4}}.
\end{align*}
\end{Theorem}

\textit{Proof.} Just use the definition of the $q$-gamma function together with the transformation formula \eqref{qq2toq} and the special values \eqref{PhiPp2}-\eqref{PhiP8}.
\hfill$\Box$

\section{Special values for the first Jacobi theta function}

Applying Theorem \ref{theta1sinq} and the special values of the $q$-sine function in Theorem \ref{sinqspecval}, we arrive at the next statement.

\begin{Theorem}\label{theta1lnq}For any $0<q<1$ we have
\begin{eqnarray}
i\vartheta_1\left(\frac12i\ln q,q\right)&=&q^{-1/4}\frac{\phi^2(q)}{\phi(q^2)},\label{spec1}\\
i\vartheta_1\left(\frac14i\ln q,q\right)&=&\frac{\phi(\sqrt{q})\phi(q^2)}{\phi(q)},\label{spec2}\\
i\vartheta_1\left(\frac18i\ln q,q\right)&=&\frac{1}{\sqrt{2}}\frac{\phi^3(q)}{\phi(q^2)\phi(\sqrt{q})}\sqrt{\frac{\phi(\sqrt[4]{q})}{\phi(\sqrt{q})}}\sqrt{\sqrt{\frac{4}{S^2(\sqrt{q})}+1}-1}.\label{spec3}
\end{eqnarray}
\end{Theorem}

After the previous steps we are ready now to list some special values for $\vt_1$. We just set $e^{-\pi/2}$, $e^{-\pi}$, $e^{-2\pi}$ and $e^{-4\pi}$ in \eqref{spec1} and calculate the given Euler function values. Then a considerable simplification can be made to get
\begin{align*}
i\vt_1\left(-\frac{i\pi}{4},e^{-\pi/2}\right)&=\frac{e^{\pi /8}}{2^{19/24}\pi ^{3/4}} \Gamma \left(\frac{1}{4}\right)\left(\sqrt{2}-1\right)^{2/3} \sqrt[12]{4+3 \sqrt{2}},\\
i\vt_1\left(-\frac{i\pi}{2},e^{-\pi}\right)&=\frac{e^{\pi /4} \Gamma \left(\frac{1}{4}\right)}{(2 \pi )^{3/4}},\\
i\vt_1\left(-i\pi,e^{-2\pi}\right)&=\frac{e^{\pi /2} \Gamma \left(\frac{1}{4}\right)}{2^{5/8} \pi ^{3/4}},\\
i\vt_1\left(-2i\pi,e^{-4\pi}\right)&=\frac{e^{\pi } \Gamma \left(\frac{1}{4}\right)}{2^{15/16}\pi ^{3/4} \sqrt[4]{\sqrt{2}-1}}.
\end{align*}

Then substituting $e^{-\pi}$, $e^{-2\pi}$, $e^{-4\pi}$ in place of $q$ in \eqref{spec2}, another set of special values arises:
\begin{align*}
i\vt_1\left(-\frac{i\pi}{4},e^{-\pi}\right)&=\frac{e^{\pi /16}}{\pi ^{3/4}} \Gamma \left(\frac{5}{4}\right)2 \sqrt[8]{3 \sqrt{2}-4},\\
i\vt_1\left(-\frac{i\pi}{2},e^{-2\pi}\right)&=e^{\pi /8} \left(\frac{2}{\pi }\right)^{3/4} \Gamma \left(\frac{5}{4}\right),\\
i\vt_1\left(-i\pi,e^{-4\pi}\right)&= \frac{e^{\pi /4}2^{9/16}}{\pi ^{3/4}} \Gamma \left(\frac{5}{4}\right) \sqrt[4]{\sqrt{2}-1}.
\end{align*}

We exhaust all of our possibilities if we substitute $q=e^{-2\pi}$ and $q=e^{-4\pi}$ into \eqref{spec3}. First we calculate $C(p)$:
\begin{align*}
C(e^{-\pi/2})&=\frac{2^{5/12}}{(\sqrt 2-1)^{2/3}(4+3\sqrt2)^{1/12}},\\
C(e^{-\pi})&=2^{5/4},\\
C(e^{-2\pi})&=\frac{2}{\sqrt 2-1}.\\
\end{align*}
This results
\begin{align*}
i\vt_1\left(-\frac{i\pi}{4},e^{-2\pi}\right)&=2^{5/8}\frac{e^{\pi /32}}{\pi ^{3/4}}\Gamma \left(\frac{5}{4}\right) \sqrt[16]{-68+323 \sqrt{2}-56 \sqrt{68-14 \sqrt{2}}},\\
i\vt_1\left(-\frac{i\pi}{2},e^{-4\pi}\right)&=\frac{e^{\pi /16}}{\pi ^{3/4}}\Gamma \left(\frac{5}{4}\right)\sqrt[4]{2+6 \sqrt{2}-4 \sqrt{2 \left(2+\sqrt{2}\right)}} .
\end{align*}

The restricted set of known values of the Euler function permits us to calculate just the above special values, not others. Employing newer techniques, however, it is possible to proceed: see the last section.

We have to note that behind the above calculations there is an unproven formula, \eqref{qcosineaddM}. However, the above special values seems to be entirely correct. As a justification, we calculated $\vt_1\left(-\frac{i\pi}{2},e^{-4\pi}\right)$ and the given expression on the right up to 500 decimal digits. There is a full agreement.

\section{A combinatorial identity between two classes of elliptic series}

S. Ramanujan investigated sums of the form
\[S_n=\sum_{m=0}^\infty(-1)^mq^{m(m+1)/2}(2m+1)^n,\]
see \cite[p.61-]{RNB3}. He also introduced the quantity
\[Q_n(q)=\frac{1}{n+1}\frac{\sum_{m=0}^\infty(-1)^mq^{m(m+1)/2}(2m+1)^{n+1}}{\sum_{m=0}^\infty(-1)^mq^{m(m+1)/2}(2m+1)}.\]
It is proven (see also \cite{Ram1,Ram2}) that $Q_n$ is always a polynomial of three functions, $L,M,N$, where
\begin{align*}
L(q)&=1-24\sum_{n=1}^\infty\frac{nq^n}{1-q^n},\\
M(q)&=1+240\sum_{n=1}^\infty\frac{n^3q^n}{1-q^n},\\
N(q)&=1-504\sum_{n=1}^\infty\frac{n^5q^n}{1-q^n}.
\end{align*}
For example \cite[p. 65]{RNB3},
\begin{align}
Q_2&=\frac13L,\\
Q_4&=\frac{1}{15}(5L^2-2M).\nonumber
\end{align}

It is a well known fact that the above sums defining $L,M,N$ are equal to
\[\sum_{n=1}^\infty\frac{n^aq^n}{1-q^n}=\sum_{n=1}^\infty\sigma_a(n)q^n,\]
where
\[\sigma_a(n)=\sum_{d|n}d^a\]
is the generalized divisor function. More on these sums can be found in \cite{Apostol}.

We introduce the sums
\begin{align*}
A(n,q)&:=\frac{1}{\vt_2(q)}\sum_{m=0}^\infty(-1)^mq^{(m+1/2)^2}(2m+1)^n,
\intertext{and}
B(n,q)&:=\frac{1}{\vt_2(q)}\sum_{m=0}^\infty q^{(m+1/2)^2}(2m+1)^n.
\end{align*}

Note that $A(n,q)$ equals to $S_n$ up to the factor $\vt_2(q)$. Hence one can rewrite $Q_n(q)$ as
\[Q_n(q)=\frac{1}{n+1}\frac{A(n+1,\sqrt{q})}{A(1,\sqrt{q})}.\]

We prove a combinatorial sum which connects the elliptic sums $A$ and $B$ which will permit us to calculate several quotients of theta functions. This sum is as follows.

\begin{Theorem}\label{combid}For all $n=0,1,2,\dots$ and for any $0<q<1$
\[C(p)\sum_{k=0}^n\binom{2n}{2k}\frac{1}{k+1}A(2k+1,\sqrt{q})B(2n-2k,\sqrt{q})=\frac{4^n}{2n+1}A(2n+1,q),\]
where $p$ and $q$ are connected via the transcendental equation
\[\ln p\ln q=\pi^2,\]
and $C(p)$ is defined in \eqref{cdef}.
\end{Theorem}

\textit{Proof.} The proof is easy, we just apply the definition of the $q$-trigonometric functions and the addition theorem \eqref{qsineaddM}. By \eqref{qsinth} and by the representations \eqref{th1ssum}-\eqref{th2ssum} we have that
\begin{align*}
\sin_p(z)&=2\sum_{n=0}^\infty\frac{(-1)^n}{(2n+1)!}z^{2n+1}A(2n+1,q),\\
\cos_p(z)&=2\sum_{n=0}^\infty\frac{(-1)^n}{(2n)!}z^{2n}B(2n,q).\\
\end{align*}
We also applied the standard Taylor series representations of  $\sin$ and $\cos$. If we look at the Taylor coefficients of
\[\sin_p(2z)\quad\mbox{and}\quad C_p\sin_{p^2}(z)\cos_{p^2}(z),\]
they must be equal by \eqref{qsineaddM}. This equality yields our theorem.
\hfill$\Box$

Note that the above sum in the theorem is the ``elliptic analogue'' of the well known finite combinatorial identity
\[\sum_{k=0}^n\binom{2n}{2k}\frac{1}{2k+1}=\frac{4^n}{2n+1}\quad(n\ge0).\]
The validity of this identity can be proven as we did above using the classical trigonometric functions instead.

In the following section we investigate the case $n=0$ in more detail.

\section{Special values of quotients of the second Jacobi elliptic function}

Theorem \ref{combid} specializes to
\[C(p)A(1,\sqrt{q})B(0,\sqrt{q})=A(1,q),\]
when $n=0$. $A(1,q)$ can be expressed with the Euler $\phi$ function:
\[A(1,q)=q^{1/4}\vt_2(q)\phi^3(q^2),\]
see \cite[p. 53]{RNB3}. Hence we win the next relation.

\begin{Theorem}\label{thetafrac}If $0<p,q<1$ such that $\ln p\ln q=\pi^2$, then
\[\frac{\vt_2(\sqrt{q})}{\vt_2(q)}=q^{-1/8}\frac{C(p)}{2}\frac{\phi^3(q)}{\phi^3(q^2)}.\]
\end{Theorem}

\begin{Corollary}If we set $q=e^{-\pi/2}$, $q=e^{-\pi}$, $q=e^{-2\pi}$ in the previous theorem, we get
\begin{align*}
\frac{\vt_2(e^{-\pi/4})}{\vt_2(e^{-\pi/2})}&=\sqrt[8]{8+6\sqrt{2}},\\
\frac{\vt_2(e^{-\pi/2})}{\vt_2(e^{-\pi})}&=2^{5/8},\\
\frac{\vt_2(e^{-\pi})}{\vt_2(e^{-2\pi})}&=\frac{2^{13/24}}{(\sqrt{2}-1)^{2/3}(4+3\sqrt{2})^{1/12}}.
\end{align*}
\end{Corollary}

This is all what we get employing the special values of $\phi$. However, we can go much more forward if we realize the connection of the right hand side of the theorem with some new theta identites. But before, we give a connection with the Weber functions.

\subsection{The Weber function $\mathfrak{f}_1$}

Theorem \eqref{thetafrac} easily yields that there is a relationship between the theta quotient $\frac{\vt_1(\sqrt q)}{\vt_1(q)}$ and the $\mathfrak{f}_1$ Weber function. The three Weber functions are defined as \cite{BB,Hart1,Hart2,Hart3,Weber}
\[\mathfrak{f}(\tau)=e^{-i\pi/24}\frac{\eta\left(\frac{\tau+1}{2}\right)}{\eta(\tau)},\quad\mathfrak{f}_1(\tau)=\frac{\eta\left(\frac{\tau}{2}\right)}{\eta(\tau)},\quad,\mathfrak{f}_2(\tau)=\sqrt2\frac{\eta(2\tau)}{\eta(\tau)}.\]
Here
\[\eta(\tau)=q^{1/24}\phi(q),\quad (q=e^{2\pi i\tau})\]
is the Dedekind eta function.

Theorem \eqref{thetafrac} readily gives that
\[\frac{\vt_2(\sqrt{q})}{\vt_2(q)}=\frac{C(p)}{2}\mathfrak{f}_1(2\tau)\quad(q=e^{2\pi i\tau}).\]
If we look at the definition of $C(p)$ and we employ the definition of $\mathfrak{f}_1$ again, we get the next theorem.

\begin{Theorem}Let $0<p,q<1$ such that $\ln p\ln q=\pi^2$. Then
\[\frac{\vt_2(\sqrt{q})}{\vt_2(q)}=\frac12\frac{\mathfrak{f}_1^4(4\tau_p)}{\mathfrak{f}_1^2(2\tau_p)}\mathfrak{f}_1^3(2\tau_q),\]
where $\tau_p$ and $\tau_q$ are defined by the equations
\[q=e^{2\pi i\tau_q},\quad p=e^{2\pi i\tau_p}.\]
\end{Theorem}

\section{Some additional special values obtained from modular equations}

As we said before, it is possible to extend our results by using Weber functions and some new theta identities obtained by several recent authors \cite{Naika,VS}.

Up to now we deduced all the possible special values what the Euler function values \eqref{PhiPp2}-\eqref{PhiP8} enabled (see the results after Theorem \ref{theta1lnq} and Theorem \ref{thetafrac}). From now on we just can give some directions what one could follow to obtain new special values for theta functions.

\subsection{The special value of $\phi(e^{-\pi/2})$}\label{Hartres}

First we show how one can proceed winning new special values of the Euler modular function $\phi$. Especially, we show that \eqref{PhiPp2} really holds:
\[\phi\left(e^{-\pi/2}\right)=\frac{\sqrt[3]{\sqrt{2}-1} \sqrt[24]{4+3 \sqrt{2}} e^{\pi /48} \Gamma \left(\frac{1}{4}\right)}{2^{5/6} \pi ^{3/4}}.\]

The Weber functions are related as follows \cite[p. 69]{BB}
\[\mathfrak{f}\mathfrak{f}_1\mathfrak{f}_2=\sqrt{2},\]
\[\mathfrak{f}^8=\mathfrak{f}_1^8+\mathfrak{f}_2^8.\]

Then substituting into this latter equation the identity $\mathfrak{f}=\frac{\sqrt 2}{\mathfrak{f}_1\mathfrak{f}_2}$ we get that
\[\left(\frac{\sqrt2}{\mathfrak{f}_1\mathfrak{f}_2}\right)^8=\mathfrak{f}_1^8+\mathfrak{f}_2^8.\]
Then employing the definition of $\mathfrak{f}_2$ we substitute
\[\mathfrak{f}_2=\frac{\sqrt2}{\mathfrak{f}_1(2\tau)}\]
and we rearrange to arrive at
\begin{equation}
\mathfrak{f}_1^{16}(2\tau)=\mathfrak{f}_1^8(2\tau)\mathfrak{f}_1^{16}(\tau)+16\mathfrak{f}_1^8(\tau).\label{f1eq}
\end{equation}
Now let $\tau=i/2$. By the relation between $\phi$, $\eta$ and $\mathfrak{f}_1$ we get that $\mathfrak{f}_1(i)=2^{1/8}$. Substituting this into \eqref{f1eq}, we can solve the resulting equation of degree 16 by \texttt{Mathematica}, for example. This gives that
\[\mathfrak{f}_1\left(\frac{i}{2}\right)=(3\sqrt 2-4)^{1/8}.\]
Hence we can trace back this by the definition of $\mathfrak{f}_1$ and the relation between $\eta$ and $\phi$ to finalize the proof. (During the course we need the special values $\eta(i)$ and $\eta(i/2)$ but these can easily be calculated from $\phi(e^{-2\pi})$ and $\phi(e^{-\pi})$, respectively.)

Note that this argument above is capable to give $\phi(e^{-\pi/4})$ etc., too.

\subsection{A modular equation of Naika}

We shall prove two additional results, namely,
\begin{equation}
\frac{\vt_2(e^{-\pi/3})}{\vt_2(e^{-2\pi/3})}=\label{thetaq}
\end{equation}
\[\frac{\left(594100+420099 \sqrt{2}\right)^{1/12} \left(1+\sqrt{3}+\sqrt{2} 3^{3/4}\right)^{4/3} \left(-2-\sqrt{3}+\sqrt{9+6 \sqrt{3}}\right)}{2^{13/24} \left(-6+5 \sqrt{2}+6\ 3^{1/4}+9 \sqrt{2} 3^{1/4}-4\ 3^{3/4}+\sqrt{2} 3^{3/4}+7 \sqrt{6}\right)^{2/3}},\]
and
\begin{equation}
\frac{\vt_1(-3i\pi/2,e^{-3\pi})}{\vt_1(-3i\pi/4,e^{-3\pi})}=\label{thetaq2}
\end{equation}
\[\frac{\left(1188200+840198 \sqrt{2}\right)^{1/24} \left(1+\sqrt{3}+\sqrt{2} 3^{3/4}\right)^{2/3} e^{9 \pi /16}}{\left(-6+5 \sqrt{2}+6\ 3^{1/4}+9 \sqrt{2} 3^{1/4}-4\ 3^{3/4}+\sqrt{2} 3^{3/4}+7 \sqrt{6}\right)^{1/3}}.\]

Following M. S. M. Naika, we define the functions
\[P=\frac{\phi(q)}{q^{1/24}\phi(q^2)},\quad Q_n=\frac{\phi(q^n)}{q^{n/24}\phi(q^{2n})}.\]
Then \cite{Naika}
\[(PQ_3)^3+\frac{8}{(PQ_3)^3}=\left(\frac{Q_3}{P}\right)^6-\left(\frac{P}{Q_3}\right)^6.\]

Putting $q=e^{-\pi/2}$ it is possible to calculate $P$, and from the equation of Naika, $Q_3$ comes by using some heavier algebra which was aided by the mathematical package \texttt{Mathematica}. Substituting the calculated value of $Q_3$ into the definition of $C(p)$ (when $\ln p\ln q=\pi^2$), \eqref{thetaq} can be calculated after some algebra by Theorem \ref{thetafrac}. The same can be done to get the special value \eqref{thetaq2} but here we have to use the quotient of \eqref{spec1} and \eqref{spec2}.

We note that Naika deduced modular equations with respect to $Q_5$, $Q_7$ and $Q_{11}$, and Vasuki and Sreeramamurthy for other variants, too. With these and with \eqref{mid} in principle is possible to get closed forms for other theta quotients, however, the underlying algebraic manipulations are getting more and more difficult.

\begin{center}
\textbf{Acknowledgement}
\end{center}

The author is grateful to Prof. William B. Hart, who gave the argument to deduce \eqref{PhiPp2} in subsection \ref{Hartres}.

\end{document}